\documentclass[10pt, twoside, a4paper]{article}
\let\finishall\relax\let\Finishall\relax\let\getprepared\relax
\ifx\TestIngCommand\undefined\relax
\else\input ../../../proc/ppreamb-book.tex 
  \input init.tex \fi
\let\TestIngCommand\undefined

\newtheorem{example}{Example}

\newtheorem{question}{Question}

\newtheorem{conj}{Conjecture}

\usepackage{amsmath,amscd,amssymb}                                         
\usepackage{graphics}                                                      
\newtheorem{theo}{Theorem}                                                 
\newtheorem{lem}{Lemma}                                                    
\newtheorem{cor}{Corollary}                                                
\newtheorem{defi}{Definition}                                              
\newskip\ttglue\ttglue=.5em plus.25em minus.15em                           
\def\firstname#1{\def\FIRSTNAME{#1}\ignorespaces}
\def\lastname#1{\def\LASTNAME{#1}\ignorespaces}
\def\middleinitial#1{\def\MIDDLEINI{#1}\ignorespaces}
\def\department#1{\def\DEPARTMENT{#1}\ignorespaces}
\def\institute#1{\def\INSTITUTE{#1}\ignorespaces}
\def\address#1{\def\ADDRESS{#1}\ignorespaces}
\def\country#1{\def\COUNTRY{#1}\ignorespaces}
\def\otheraffiliation#1{\def\OTHERAFFILIATION{#1}\ignorespaces}
\def\email#1{\def\EMAIL{#1}\ignorespaces}
\newcount\autcount\autcount=0                                              
\newcount\affcount\affcount=0                                              
\newcount\numcount\newcount\nummcount                                      
\newcount\nummmcount\newcount\nummmmcount                                  
\def\writename#1#2{\ \kern-1ex\hbox{
  \csname AUthor\the#1\endcsname\                                          
  \edef\TESTSTR{}\expandafter\ifx\csname auTHor\the#1\endcsname\TESTSTR    
  \else\csname auTHor\the#1\endcsname.\ \fi                                
  \csname authOR\the#1\endcsname$^{\csname AFF\the#1\endcsname}$
  \expandafter\ifx\csname corr\number#1\endcsname\relax                    
  \else\thanks{Corresponding author.}\ \fi                                 
  }\ifnum#1<#2, \else\ \kern-1ex\fi}
\def\writeemail#1{
  \nummcount=0\relax\nummmcount=0\relax                                    
  \loop\ifnum\nummcount<\autcount\advance\nummcount by1\relax              
    {\expandafter\ifnum\csname AFF\the\nummcount\endcsname=#1\relax        
    \global\advance\nummmcount by1\fi}\repeat                              
  \nummcount=0\relax\nummmmcount=0\relax                                   
  \loop\ifnum\nummcount<\autcount\advance\nummcount by1\relax              
    {\expandafter\ifnum\csname AFF\the\nummcount\endcsname=#1\relax        
    \global\advance\nummmmcount by1\relax\def\blank{}\expandafter          
    \ifx\csname EMAIL\the\nummcount\endcsname\blank(no e-mail)
    \else\csname EMAIL\the\nummcount\endcsname                             
    \fi                                                                    
    \ifnum\nummmmcount<\nummmcount; \fi\fi}\repeat}
\long\def\BeginAuthorList#1\EndAuthorList{#1\relax                         
  \author{\vbox{\hsize=390pt\noindent\numcount=0\relax                     
    \loop\ifnum\numcount<\autcount\advance\numcount by1\relax              
      \writename{\numcount}{\autcount}
      \repeat}\\[2mm]                                                      
    \vbox{\small\numcount=0\relax                                          
      \loop\ifnum\numcount<\affcount\advance\numcount by1\relax            
        \vbox{{\count0=\numcount\relax                                     
          \loop\expandafter\ifnum\csname AFF\the\count0\endcsname
            <\numcount\relax\advance\count0 by1\relax\repeat               
          $^{\csname AFF\the\count0\endcsname}$}
        \def\BLANK{}\expandafter\ifx\csname DEPT\the\numcount\endcsname    
          \BLANK                                                           
          \else\csname DEPT\the\numcount\endcsname, \fi                    
        \csname INST\the\numcount\endcsname,                               
        \csname ADDR\the\numcount\endcsname,                               
        \csname COUN\the\numcount\endcsname                                
        \edef\TEST{}\expandafter\ifx\csname OTHE\the\numcount\endcsname
          \TEST                                                            
          .\else;\break\csname OTHE\the\numcount\endcsname.\fi}
        \vbox{\writeemail{\numcount}}
        \repeat}\\}}
\expandafter\def\csname x1\endcsname{}
\expandafter\def\csname x2\endcsname{}
\expandafter\def\csname x3\endcsname{}
\expandafter\def\csname x4\endcsname{}
\expandafter\def\csname x5\endcsname{}
\expandafter\def\csname x6\endcsname{}
\expandafter\def\csname x7\endcsname{}
\expandafter\def\csname x8\endcsname{}
\expandafter\def\csname x9\endcsname{}
\def\Author#1#2{\global\advance\autcount by1\relax#2                       
  \expandafter\edef\csname AUthor\the\autcount\endcsname{\FIRSTNAME}
  \expandafter\edef\csname auTHor\the\autcount\endcsname{\MIDDLEINI}
  \expandafter\edef\csname authOR\the\autcount\endcsname{\LASTNAME}
  \expandafter\edef\csname EMAIL\the\autcount\endcsname{\EMAIL}
  \let\tempera\"\def\"{\string\"}\expandafter\ifx\csname x\DEPARTMENT
    \endcsname\relax                                                       
    \global\advance\affcount by1\relax\let\"\tempera                       
    \expandafter\edef\csname DEPT\the\affcount\endcsname{\DEPARTMENT}
    \expandafter\edef\csname INST\the\affcount\endcsname{\INSTITUTE}
    \expandafter\edef\csname ADDR\the\affcount\endcsname{\ADDRESS}
    \expandafter\edef\csname COUN\the\affcount\endcsname{\COUNTRY}
    \expandafter\edef\csname OTHE\the\affcount\endcsname{\OTHERAFFILIATION}
    \expandafter\edef\csname AFF\the\autcount\endcsname{\the\affcount}
  \else\expandafter\edef\csname AFF\the\autcount\endcsname{\DEPARTMENT}
  \fi\let\"\tempera\ignorespaces}
\def\CorrespondingAuthor#1#2{
  \expandafter\xdef\csname corr\number#1\endcsname{cor}
  \Author#1{#2}}
\def\PaperTitle#1{\title{\bf#1}}
\def\Category#1{\ignorespaces}
\def\keywords#1{{\noindent \emph{Keywords:}                                
  \def\BLANK{}\def\TEST{#1}\ifx\BLANK\TEST(n/a).\else#1\fi}}
\getprepared                                                               
\begin{document}                                                           
\PaperTitle{Bose and Einstein Meet Newton}%
\Category{(Pure) Mathematics}

\date{}
\BeginAuthorList 
  \Author1{
    \firstname{Wayne}
    \lastname{Lawton}
    \middleinitial{M}   
    \department{Department of Mathematics}
    \institute{Mahidol University}
    \address{Bangkok 10400}
    \country{Thailand}
    \otheraffiliation{School of Mathematics and Statistics, University of Western Australia, Perth, Australia}
    \email{wayne.lawton@uwa.edu.au}}
\EndAuthorList 
\maketitle 
\thispagestyle{empty} 
\begin{abstract} 
We model the time evolution of a Bose-Einstein condensate, subject to a special periodically excited optical lattice,
by a unitary quantum operator $U$ on a Hilbert space $H.$ If a certain parameter $\alpha = p/q,$ where $p$ and $q$ are coprime positive integers, then
$H = L^2(\mathbb{R}/\mathbb{Z},\mathbb{C}^q)$
and $U$ is represented by a $q \times q$ matrix-valued function $M$ on
$\mathbb{R}/\mathbb{Z}$
that acts pointwise on functions in $H.$ The dynamics of the quantum system is described by the eigenvalues of $M.$ Numerical computations
show that the characteristic polynomial $\det(zI - M(t)) = \prod_{j=1}^{q} \left(z - \lambda_j(t) \right)$ where each
$\lambda_j$ is a real analytic function that has period $1/q.$ We discuss this phenomena using Newton's Theorem, published in {\it Geometria analytica} in 1660, and modern concepts from analytic geometry.
\end{abstract} 
\keywords{bose-einstein condensate, unitary quantum operator, characteristic polynomial,
newton's theorem, resolution of singularities, \'{e}tale homotopy.}
\finishall 
\section{Introduction and Preliminary Results}
This paper discusses polynomials, which are used to model
the dynamics of certain quantum systems involving Bose-Einstein condensates
controlled by optical lattices, and that arise as follows.
\\ \\
Let $\alpha = p/q$ where $p = 2$ and $q = 5,$
and construct the family of matrix valued functions
$M_\kappa(t) = D_\kappa(t) \, G^{-1} \, D_\kappa(t) \, G, \ \ \kappa > 0$ where
$$
G = \left[
  \begin{array}{ccccc}
    g_0 & g_1 & g_2 & g_3 & g_4 \\
    g_4 & g_0 & g_1 & g_2 & g_3 \\
    g_3 & g_4 & g_0 & g_1 & g_2 \\
    g_2 & g_3 & g_4 & g_0 & g_1 \\
    g_1 & g_2 & g_3 & g_4 & g_0 \\
  \end{array}
\right], \ D_{\kappa}(t) = \left[
  \begin{array}{ccccc}
    c_0(t)  & 0 & 0 & 0 & 0 \\
    0 & c_1(t) & 0 & 0 & 0 \\
    0 & 0 & c_2(t) & 0 & 0 \\
    0 & 0 & 0 & c_3(t) & 0 \\
    0 & 0 & 0 & 0 & c_4(t) \\
  \end{array}
\right],
$$
$g_j = \frac{1}{q} \, \sum_{k = 0}^{q-1} e^{-i2\pi k^2 \alpha} \, e^{i2\pi jk \alpha}$
and $c_j(t) = e^{-i 2 \kappa \cos 2\pi (t-j\alpha)}$ for $j = 0, 1, 2, 3, 4.$
\\ \\
As explained in Appendix A, the matrix valued function $M_\kappa$ represents a Floquet operator that describes the time evolution of a
quantum system called an {\it on resonance double kicked rotor}.
The characteristic polynomial $C(t,z) = \det(zI - M_\kappa(t))$ can
be regarded as a polynomial of degree 5 in the $z$ variable whose coefficients are functions of the $t$ variable that are analytic and
have period $1.$ A closer inspection of the structure of $M_\kappa$ reveals that each coefficient of $C(t,z)$ also has period $\alpha = p/q$ and
therefore has period $1/q.$ Therefore there exists a polynomial $P(t,z)$ of $z$ whose coefficients are functions of $t$ that have
period $1$ such that $P(qt,z) = C(t,z).$ The dynamics of the system is described by the spectrum of the
Floquet operator which equals $\cup_{t=0}^{1} \hbox{roots} P(t,z).$ This implies that the spectrum consists of a union
of at most $q = 5$ disjoint closed intervals or {\it bands}. Numerical investigations show that for sufficiently small
values of $\kappa$ there are $q$ bands for odd $q$ and $q-1$ bands for even $q.$ This can be shown to imply that each root
of $P(t,z)$ is a continuous period $1$ function of $t.$
As $\kappa$ increases the crossing of the graphs of the roots increases and for sufficiently large $\kappa$ there is only $1$ band.
However, if we fix the value of $\kappa$ and construct matrices corresponding to rational $\alpha = p/q$ as $\alpha$ approaches an
irrational number, and therefore $q \rightarrow \infty,$ we find that more bands appear. Based on extensive numerical experiments
we conjectured in a previous paper \cite{lmwg} that the spectrum approaches a Cantor set. Strong support for our Cantor conjecture would
be provided by
\begin{conj}
\label{Conj}
If $\alpha$ is rational and $\kappa > 0$ and $P$ is constructed as above for the Floquet operator corresponding to $\alpha$ and $\kappa,$ then each root of $P(t,z)$ is an analytic function of $t$ that has period $1.$
\end{conj}
This conjecture is supported by extensive numerical computation. The results in this paper imply by a homotopy argument that if the roots of $P(t,z)$ are locally analytic functions then the conjecture holds. They also imply that if the multiplicity of the roots is at most 2 then the conjecture holds. A proof of this conjecture would support our Cantor conjecture since if each root of $P(t,z)$ is an analytic function of period $1$ then each root of $C(t,z)$ is analytic and has period $1/q$ so this gives a small upper bound on the range of each root.
Since this is a pure mathematics paper we summarize the physical background in Sections 5 and 6 where references
for further study are suggested.
\section{Preliminary Results}
$\mathbb{Z}, \mathbb{Q}, \mathbb{R}, \mathbb{C}$ denote the integer, rational,
real, and complex numbers, $\mathbb{T} = \mathbb{R}/\mathbb{Z}$ denotes the real circle
group.
For any algebra $A,$ $A[z]$ denotes the algebra of polynomials with coefficients in $A.$
A polynomial is called {\it monic} if the coefficient of $z^{\hbox{deg}(P)}$ equals $1,$
{\it simple} if there does not exist $Q \in A[z]$ with deg$(Q) \geq 1$ such that $Q^2$ divides $P,$ {\it irreducible} if $P$ does not admit a factorization $P = QR$ where deg$(Q) \geq 1$ and deg$(R) \geq 1,$ and {\it completely reducible} (CR) if it admits a factorization
$P(z) = \prod_{j=1}^{\deg P} (z-a_j)$ with $a_j \in A.$
$K = C^{\, \omega}(\mathbb{T})$ denotes the algebra of analytic
functions on $\mathbb{T}$ and $K_0 = C_{0}^{\omega}(\mathbb{T})$ denotes the algebra of {\it germs of analytic functions} at $0$ (power series in $t$ that converge absolutely for $t$ sufficiently small). Clearly $K \subset K_0$ and $K \neq K_0.$ For $P \in K[z],$ $P_0 \in K_{0}[z]$ denotes the element obtained by regarding the coefficients of $P$ to be in $K_0.$ %
$\mathbb{T}$ acts as a group of algebra automorphisms of $K[z]$ by rotation
$
    (R_s \, P)(t,z) = P(t+s,z).
$
For $P \in K[z]$ and $s \in \mathbb{T}$ we define $P_s = (R_sP)_{0} \in K_{0}[z].$
We call $P \in K[z]$ {\it locally completely reducible} (LCR) if for every $s \in \mathbb{T},$ $(R_sP)_0 \in K_0[z]$ is CR.

\begin{example}
\label{Example 1} If $P(t,z) = z^2 - \sin^2 2\pi t$ then $P(0,z) \in \mathbb{C}[z]$ is not simple. However $P \in K[z]$ and $P_0 \in K_0[z]$ are simple.
\end{example}

\begin{question}
\label{Question 1}
When is a monic polynomial in $A[z]$ simple for $A = \mathbb{C}, \, A = K, \, A = K_0 \, ?$
\end{question}

\begin{example}
\label{Example 2}
$P(t,z) = z^2 - \cos 2\pi t \in K[z]$ is not LCR since $P_{\pm 1/2} \in K_0[z]$ are not CR.
\end{example}

\begin{question}
\label{Question 2}
When is a monic polynomial $P \in K[z]$ LCR?
\end{question}

\begin{example}
\label{Example 3}
If $P(t,z) = z^2 - e^{2\pi i t} \in K[z]$ then $P$ is LCR but not CR.
\end{example}

\begin{question}
\label{Question 3}
When is a LCR monic polynomial $P \in K[z]$ CR?
\end{question}
If $A$ is one of the algebras $\mathbb{C},$ $K,$ or $K_0,$ then $A[z]$ is a {\it Euclidean domain}, therefore $P(z) = p_mz^m + p_{m-1}z^{m-1} + \cdots + p_1z+p_0$ and  $Q(z) = q_nz^n + q_{n-1}z^{n-1} + \cdots + q_1z + q_0$ in $A[z]$ have a {\it greatest common divisor}
GCD$(P,Q).$ Furthermore, the {\it Euclidean Algorithm}, codified about 300 BCE by
Euclid in Books VII and X in his Elements \cite{heath}, but likely known to Eudoxus of Cnidus about 375 BCE \cite{becker}, for computing the greatest common divisor of two positive integers, can be also used to compute $P_1, Q_1 \in A[z]$ such that GCD$(P,Q) = P_1P + Q_1Q$ (perhaps the first Bezout identity \cite{bezout}). This gives an algorithmic solution to Question 1. We now give an explicit solution. The {\it Sylvester Resultant} $R(P,Q)$ is the determinant of the following $(m+n) \times (m+n)$ matrix
$$
\left[
\begin{array}{ccccccccccccc}
p_m     & p_{m-1}   & \cdot     & \cdot     & \cdot     & \cdot     & \cdots & \cdot     & p_1       & p_0       & 0         & \cdots    & 0 \\
0       & \ddots    & \ddots    & \ddots    & \ddots    & \ddots    & \ddots & \ddots    & \ddots    & \ddots    & \ddots    & \ddots    & \vdots \\
\vdots  & \ddots    & \ddots    & \ddots    & \ddots    & \ddots    & \ddots & \ddots    & \ddots    & \ddots    & \ddots    & \ddots    & 0 \\
0       & \cdots    &  0        & p_{m}     & p_{m-1}   & \cdot     & \cdots & \cdots    & \cdot     & \cdot     & \cdot     & p_1       & p_0 \\
q_n     & q_{n-1}   &  \cdot    & q_1       & q_0       & 0         & \cdots & \cdots    & \cdots    & \cdots    & \cdots    & \cdots    & 0 \\
0       & \ddots    & \ddots    & \ddots    & \ddots    & \ddots    & \ddots & \ddots    & \ddots    & \ddots    & \ddots    & \ddots    & \vdots \\
\vdots  & \ddots    & \ddots    & \ddots    & \ddots    & \ddots    & \ddots & \ddots    & \ddots    & \ddots    & \ddots    & \ddots    & \vdots \\
\vdots  & \ddots    & \ddots    & \ddots    & \ddots    & \ddots    & \ddots & \ddots    & \ddots    & \ddots    & \ddots    & \ddots    & \vdots \\
\vdots  & \ddots    & \ddots    & \ddots    & \ddots    & \ddots    & \ddots & \ddots    & \ddots    & \ddots    & \ddots    & \ddots    & \vdots \\
\vdots  & \ddots    & \ddots    & \ddots    & \ddots    & \ddots    & \ddots & \ddots    & \ddots    & \ddots    & \ddots    & \ddots    & \vdots \\
\vdots  & \ddots    & \ddots    & \ddots    & \ddots    & \ddots    & \ddots & \ddots    & \ddots    & \ddots    & \ddots    & \ddots    & \vdots \\
\vdots  & \ddots    & \ddots    & \ddots    & \ddots    & \ddots    & \ddots & \ddots    & \ddots    & \ddots    & \ddots    & \ddots    & 0 \\
0       & \cdots    & \cdots    & \cdots    & \cdots    & \cdots    & \cdots & 0         & q_n       & q_{n-1}   & \cdot     & q_1       & q_0
\end{array}
\right]\, .
$$
In 1840 Sylvester proved \cite{sylvester} that if $A = \mathbb{C}$ and $p_m \neq 0$ and $q_n \neq 0$ then $R(P,Q) = 0$ if and only if $P$ and $Q$ have a common root, or equivalently,
$\hbox{deg}(\hbox{GCD}(P,Q)) \geq 1.$ It follows directly that this result also holds for $A = K$ and $A = K_0.$ The {\it derivative} of $P \in A[z]$ is the polynomial
$P^{\prime}(z) = mp_mz^{m-1} + \cdots + 2p_2z + p_1 \in A[z],$
and the {\it discriminant} of a monic $P \in A[z]$ with $\hbox{deg}(P) \geq 2$ is
$D(P) = - R(P,P^{\prime}).$
For example
\begin{equation}
\label{D}
D(z^2 + p_1z + p_0) = - \det \, \left[
                                    \begin{array}{ccc}
                                      1 & p_1 & p_0 \\
                                      2 & p_1 & 0 \\
                                      0 & 2 & p_1 \\
                                    \end{array}
                                  \right] = p_1^2 - 4p_0.
\end{equation}
For monic $P \in \mathbb{C}[z],$ $D(P) = -\prod_{i \neq j} (\lambda_i - \lambda_j),$
$\lambda_i$ are roots of $P,$ (\cite{lang}, Proposition 10.5).
\begin{lem}
\label{Q1}
A monic polynomial $P \in \mathbb{C}[z]$ is not simple iff $D(P) = 0.$ A monic polynomial $Q \in K_0[z]$ is not simple iff $D(Q) = 0$ $($this means that $D(Q)(t) = 0$ for $t$ sufficiently small$\,).$ A monic polynomial $P \in K[z]$ is not simple iff $P_0 \in K_0[z]$ is not simple.
\end{lem}
{\it Proof} The first and second assertions follow since $P$ is not simple iff $P$ and $P^{\prime}$ have a common factor with degree $\geq 1.$ The third assertion follows from the facts that $D(P(t,z))$ is an analytic function of $t$ and $D(P(t,z)) = D(P_0)(t)$ so if it vanishes for
$t$ sufficiently small then it vanishes for all $t \in \mathbb{T}.$
\\ \\
For the polynomial $P = z^2 -\sin^2 2\pi t$ in Example \ref{Example 1}, since
$D(P)= 4 \sin^2 2\pi t$ vanishes at $t = 0$ but does not vanish in a neighborhood of $0,$ $P(0,z)$ is simple but $P$ and $P_0$ are not simple. Lemma \ref{Q1} gives a complete answer to Question \ref{Question 1} and reduces Questions \ref{Question 2} and \ref{Question 3} to
questions about simple monic polynomials. Furthermore, the following result provides a partial answer to Question \ref{Question 2}.
\begin{lem}
\label{implicit}
If $Q \in K_0[z]$ is monic and $Q(0,z) \in \mathbb{C}[z]$ is simple then $Q$ is CR.
\end{lem}
{\it Proof} Since each root $\mu$ of $Q(0,z)$ has multiplicity $1,$
\begin{equation}
\label{IFT}
Q^{\prime}(0,z) = \frac{\partial Q}{\partial z}(0,\mu) \neq 0,
\end{equation}
and hence the Implicit Function Theorem (for Holomorphic Functions of Several  Complex Variables) (\cite{gunning}, Chapter 1. Theorem 4) implies that there exists $\eta_\mu \in K_0$ such that $\eta_\mu(0) = \mu$ and $Q(t,\eta_\mu(t)) = 0$ for $t$ sufficiently small. Then
$Q(t,z) = \prod_{\mu} (z-\eta_\mu(t))$
where the product is over the roots of $Q(0,z)$ so $Q$ is CR.
\begin{cor}
\label{unitary}
If $P(t,z) = z^2 + p_1(t)z + p_0(t) \in K_0[z]$ has real valued coefficients and its roots have modulus $1$ then $P$ is CR.
\end{cor}
{\it Proof}  Since the roots are complex conjugates of each other $p_0 = 1$ and $|p_1(t)| \leq 2.$ Therefore the roots have the form
$$\lambda_{\pm}(t) = -\frac{1}{2} p_1 \pm \frac{1}{2} i \sqrt {4 - p_{1}^{2}(t)}.$$
If $p_{1}^{2}(0) < 4$ then the roots are distinct so Lemma \ref{implicit} implies that
$\lambda_{\pm} \in K_0.$
Otherwise $p_{1}^{2}(t)$ has a maximum value of $4$ at $t = 0$ and hence either $p_{1}^{2}(t) = 4$ for all $t$ sufficiently small or there exists $c > 0$ and positive integer $k$ such that
$p_{1}^{2}(t) = 4 - ct^{2k} + \hbox{higher order terms}.$ Therefore $\lambda_{\pm} \in K_0.$
In either case $P(z) = (z-\lambda_{+})(z-\lambda_{-})$ is CR.
\begin{defi}
\label{primary}
We call a monic polynomial $Q \in \mathbb{C}[z]$ ${\it primary}$ if $Q(z) =  (z-\mu)^m$ for some $\mu \in \mathbb{C}.$ A monic polynomial $Q \in K_0[z]$ is called {\it point primary} $(\hbox{PP}\, )$ if $Q(0,z)$ is primary. A factorization of $Q = Q_1\, \cdots \, Q_n, \, Q_j \in K_0[z]$ is called a {\it point primary factorization} $(\hbox{PPF}\,)$ if each factor $Q_j$ is PP and
$Q_1(0,z),...,Q_n(0,z) \in \mathbb{C}[z]$
have distinct roots.
\end{defi}
$Q \in \mathbb{C}[z]$ is primary iff the {\it principal ideal}
$(Q) = \{QP : P \in \mathbb{C}[z]\}$
is a {\it primary ideal}.
The factorization of $Q$ into primary factors corresponds to the {\it primary decomposition} of $(Q),$ which is a primary topic in commutative algebra
(\cite{atiyah}, Chapter 4), (\cite{eisenbud}, Chapter 3), (\cite{kunz}, Chapter VI. Section 2), (\cite{lang}, Chapter VI, Section 5). However, PPF does not correspond to primary decomposition in $K_0[z].$
For any polynomial $Q \in \mathbb{C}[z],$ let $\Lambda(Q)$ denote the set of (distinct) roots of $Q$ and for $\mu \in \Lambda(Q),$ let $m(\mu)$ denote the multiplicity of $\mu.$
\begin{theo}
\label{PPF}
$($Hensel$\, )$ Every monic $Q \in K_0[z]$ admits a PPF.
\end{theo}
{\it Proof}
Consider the primary factorization
\begin{equation}
\label{pf}
    Q(0,z) = \prod_{\mu \in \Lambda(Q(0,z))} (z - \mu)^{m(\mu)}.
\end{equation}
Since the factors $(z-\mu)^{m(\mu)}, \, \mu \in \Lambda(Q(0,z))$ are pairwise relatively prime, Hensel's lemma \cite{hensel}, proved in 1908, implies that for every
$\mu \in \Lambda(Q(0,z))$
there exists
$Q_{\mu} \in K_0[z]$
such that
$Q_\mu(0,z) = (z - \mu)^{m(\mu)}$
and
$
Q = \prod_{\mu \in \Lambda(Q(0,z))} Q_\mu.
$
This concludes the proof. Abhyankar (\cite{abhyankar2}, 90--92) gives a algebraic proof of Hensel's lemma for polynomials with coefficients in the algebra of formal power series in an arbitrary field and gives an exercise (\cite{abhyankar2}, p. 92) that implies Hensel's lemma holds in $K_0[z].$
We give an analytic proof using {\it Cauchy's Residue Formula}, which he presented to the Academy of Sciences of Turin in 1831. Let
\begin{equation}
\label{r}
 r = \frac{1}{3} \min \{ \, |\mu - \xi| \, : \mu, \xi \in \Lambda(Q(0,z)), \  \mu \neq \xi, \, \}.
\end{equation}
For each $\mu \in \Lambda(Q(0,z))$ let
$
\Omega_\mu = \{ \, z \in \mathbb{C} \,  : \, |z-\mu| < r \, \},
$
construct the following circular contour oriented counterclockwise
$
\Gamma_\mu = \{ \, \mu + r\, e^{i\theta} \, : \, \theta \in [0,2\pi) \, \},
$
choose $\delta_\mu > 0$ so that $|Q(t,z)| > 0$ whenever $t \in (-\delta_\mu,\delta_\mu)$ and $z \in \Gamma_\mu.$ Then construct
\begin{equation}
\label{I1}
    I_\mu(t,z) = \frac{1}{2 \pi i} \int_{w \in \Gamma_\mu} \frac{Q^{\prime}(t,w)}{Q(t,w)} \, \frac{1}{z-w} \, dw, \ \ t \in (s-\delta,s+\delta), \ z \notin \Omega_\mu \cup \Gamma_\mu,
\end{equation}
and
\begin{equation}
\label{Qmu}
    Q_\mu(t,z) = \prod_{\lambda \in \Lambda(Q(0,z)) \cap \Omega_\mu} (z-\lambda)^{m(\lambda)}, \ t \in (-\delta_\mu,\delta_\mu), \, z \in \hbox{complement of } \Omega_\mu \cup \Gamma_\mu.
\end{equation}
Cauchy's residue formula implies that
\begin{equation}
\label{I2}
    I_\mu(t,z) = \frac{Q_{\mu}^{'}(z)}{Q_{\mu}(z)} = \frac{d}{dz} \log Q_\mu(t,z)
\end{equation}
and $Q(t,z) = \prod Q_\mu(t,z)$ whenever $t \in \cap_{\mu \in \Lambda(Q(0,z))} (-\delta_\mu,\delta_\mu)$ and $z$ is in the complement of $\cup_{\mu \in \Lambda(Q(0,z))} \Omega_\mu \cup \Gamma_\mu.$
Since each $I_\mu(t,z)$ is an analytic function of $t$ it follows that each $Q_\mu \in K_0[z].$ This concludes the proof.
\begin{cor}
\label{corPPF}
A monic $Q \in K_0[z]$ is CR iff each of its primary factors is CR.
\end{cor}
Lemma \ref{Q1} shows that if $Q \in K_0[z]$ is not simple then we may use the Euclidean algorithm to factorize it into simple factors. This fact in combination with Theorem \ref{PPF} can be used to reduce Question \ref{Question 2} to
\begin{question}
\label{Question 4}
When is a monic simple PP polynomial $Q \in K_0[z]$ CR?
\end{question}
Lemma \ref{Q1}, in combination with the obvious fact that if $P \in K[z]$ is CR then $P$ is LCR,
can be used to reduce Question \ref{Question 3} to
\begin{question}
\label{Question 5}
When is a LCR monic simple polynomial $P \in K[z]$ CR?
\end{question}
\section{Meet Newton}
The objective of this section is to give (at least) a partial answer to Question \ref{Question 4}.
\begin{theo}
\label{newton1}
$($Newton$\,)$ If $Q \in K_0[z]$ is a monic polynomial of degree $n > 0$ then there exists a positive integer $m$ that divides $n!$ and
$\eta_1,...,\eta_m \in K_0$ and $\delta > 0$ such that
\begin{equation}
\label{factor}
    Q(t^m,z) = \prod_{j=1}^{m} \left(z-\eta_j(t)\right), \ \ t \in (-\delta,\delta).
\end{equation}
If $Q$ is irreducible then $m = n$ and the roots can be labeled so that $\eta_j(t) = \eta_m(e^{2 \pi i j / m} t).$
\end{theo}
{\it Proof} These two statements are the Supplement 1 and Supplement 2 cases of Newton's Theorem that Abhyankar proves in (\cite{abhyankar2},89-98) using Hensel's lemma and Newton's generalized binomial theorem. He also remarks: "Newton proved this theorem about 1660 \cite{newt}. It was revived by Puiseux in 1850 \cite{puis}. The relevant history can be found in G. Chrystal's {\it Textbook of algebra}, vol. 2, 396 \cite{cry}."

\begin{example}
\label{Example 4}
If $Q = z^2 + q_1z + q_0 \in K_0[z]$ then there exists $k \geq 0,$ $c \neq 0,$ and $g \in K_0$
such that $D(Q)(t) = t^k(c+tg(t))$ and
$Q(t^2,z) = (z-\eta_1(t))(z-\eta_2(t))$ where
$$\eta_j(t) = -a_1(t^2)/2 + (-1)^j\, t^k \, (c+t^2g(t^2))^{1/2}, \ \ j = 1,2.$$ Therefore $Q$ is CR iff $k$ is even and is irreducible iff $k$ is odd. Hence if
$P(z) = z^2 - \cos 2\pi t$ is the polynomial in Example \ref{Example 2}, then $P_{\pm 1/2} \in K_0[z]$
is irreducible.
\end{example}
Equation \ref{factor} gives an analytic parameterization for the {\it analytic set}
$\{(t,z) : Q(t,z) = 0\}.$ This is an example of a {\it resolution of singularities}
that has been a central story in algebraic geometry leading to Hironaka's seminal
1964 paper \cite{hironaka}. Abhyankar \cite{abhyankar4} gives a fascinating account
of this story. Theorem \ref{newton1} implies that a monic polynomial $Q \in K_0[z]$ is CR iff the Taylor series of each $\eta_j$ in Equation \ref{factor} only has terms $ct^{k}$ where $m$ divides $k.$ Unfortunately, this fact is not very useful because it does not provide an algorithm to compute the $\eta_j.$
We need more help from Newton.
\\ \\
For any $f \in K_0$ with $f \neq 0$ let $ord\, (f) \in \{0,1,2,...\}$ denote the smallest integer $\ell$ such that $f^{(\ell)}(0) \neq 0.$ If $\eta(t)$ is a {\it formal power series} in nonnegative fractional powers of $t$ and $\eta \neq 0,$ $ord(\eta)$ denotes the smallest power of $t$ in the power series expansion. Such $\eta$ do not usually belong to $K_0.$ However, every element $K_0$ can be identified with the formal power series defined by its Taylor expansion.
Examples: $ord\, (\cos 2\pi t) = 0, ord\, (\sin 2\pi t) = 1, ord\, (t^{1/2} +t - t^{4/3}) = 1/2.$
\begin{defi}
\label{polygon}
   For $Q(z) = q_mz^m + q_{m-1}z^{n-1} + \cdots + q_1z + q_0 \in K_0[z]$ with $q_m = 1$ and $q_0 \neq 0$ let
    $\mathfrak{P}(Q) = \{ \, (j,ord\, (q_j)) \, : q_j \neq 0, \,  j = 0,..,m \}.$
    The {\it Newton Polygon} $\mathfrak{N}(Q)$ of $Q$ is the convex hull in
    $\mathbb{R}^2$ of $\mathfrak{P}(Q).$ Let $\mathfrak{E}(Q)$ denote the extreme points
    of $\mathfrak{N}(Q).$ Then $\mathfrak{E}(Q) \subseteq \mathfrak{P}(Q),$
    $(0,ord(q_0)) \in \mathfrak{E}(Q),$ and $(m,ord(q_m)) = (m,0) \in \mathfrak{E}(Q).$
    Starting with $(0,ord(q_0))$ we traverse the points in
    $\mathfrak{E}(Q)$
    in a counterclockwise direction until we reach $(m,0)$ to obtain $k+1$ points
    $(x_0,y_0) = (0,ord(q_0)),(x_1,y_1),...,(x_k,y_k) = (m,0)$ where $k \leq m.$
    This gives positive integers $m_j = x_j-x_{j-1}, \ j = 1,...,k$ which satisfy
    $m_1 + \cdots + m_k = m$ and
    slopes $s_j = (y_j-y_{j-1})/m_j, \ j = 1,...,k$ which satisfy
    $s_1 < s_2 < \cdots < s_k \leq 0.$
\end{defi}
\begin{theo}
\label{newton2}
$($Newton$\, )$ For $j = 1,...,k,$ $Q(t,z)$ has formal power series roots $($with possible multiplicity $ > 1)$ $\eta_\ell(t), \ell = 1,...,m_j$ that satisfy $ord(\eta_\ell) = -s_j.$
\end{theo}
{\it Proof} This well know result was derived by Newton in his {\it Methodus Fluxionum et
Serierum infinitarum} between 1664 et 1671 and translated into English by John Colson in 1736. See Chrystal's historical note in (\cite{cry}, Part II, p. 396) and Harold Edward's essay \cite{edwards}.
The second assertion in the following result is analogous to the PPF in Theorem \ref{PPF}.
\begin{cor}
\label{straight}
If $Q \in K_0[z]$ is irreducible then in Theorem \ref{newton2}, $k = 1,$ $m_1 = m,$
and there exists $\eta \in K_0$ with $ord(\eta) = y_1-y_0 = ord(q_0)$ such that
$\eta_\ell(t^m) = \eta(e^{2\pi i \ell/m}t).$ Furthermore,
$\prod_{\ell=1}^{m_j} (z-\eta_\ell(t)) \in K_0[z].$
\end{cor}
{\it Proof} The first assertion follows directly from Theorems \ref{newton1} and \ref{newton2}.
The second assertion follows from the fact that if $Q$ is factored into irreducible factors in
$K_0[z],$ then in Theorem \ref{newton2} for $j = 1,...,k,$ $m_j$ is the sum over the irreducible factors $P$ of $Q$ of the number of formal power series roots $\xi$ of $P$ such that $ord(\xi) = s_j.$
\\ \\
The converse of the first assertion in Corollary \ref{straight} was disproved by Abhyankar (\cite{abhyankar2}, 185--186) by constructing a reducible polynomial whose Newton polygon is a straight line. Anhyankar (\cite{abhyankar2}, p. 185) gives necessary and sufficient criteria for irreducibility in $K_0[z]$ and give a comprehensive treatment of the question of irreducibility in
\cite{abhyankar3}. The discussions in this section give (at least) a partial answer to
Question \ref{Question 4} and hence to Question \ref{Question 2}.
\section{From Jets to Braids}
The objective of this section is to completely answer Question \ref{Question 5} and thus Question \ref{Question 2}.
\begin{defi}
\label{jetdef}
For every integer $k \geq 0$ we define the k-{\it jet function}
$J_{k} : K_{0} \rightarrow \mathbb{C}^{k+1}$ by
\begin{equation}
\label{jeteq}
    J_{k}(f) = (f(0),f^{\prime}(0),f^{(2)}(0),...,f^{(k)}(0)), \ \ f \in K_0.
\end{equation}
\end{defi}
\begin{lem}
\label{jetsimple1}
If $P \in K_0[z]$ is monic, simple, CR, and $\hbox{deg}(P) = n$ then there exists
an integer $k \geq 0$ such that $J_k(\lambda_1),...,J_k(\lambda_n)$ are distinct where
$\lambda_1,...,\lambda_n \in K_0$ are the roots of $P.$
\end{lem}
{\it Proof} Assume to the contrary that for every integer $k \geq 0,$  $J_k(\lambda_1),...,J_k(\lambda_n)$ are not distinct. Then there exists $1 \leq i < j \leq n$ such that $\lambda_{i}$ and $\lambda_j$ and all of their derivatives have the same value at $0.$ Since the roots are analytic $\lambda_i = \lambda_j$ so $(z-\lambda_i)^2$ divides $P$ and hence $P$ is not simple contrary to our assumption. This contradiction concludes the proof.
\begin{example}
\label{Example 5}
$P(t,z) = z^2 - e^{2\pi i t}(1+e^{2\pi i t})z + e^{6\pi i t} \in K_0[z]$ is
simple and CR but $P(0,z) = (z-1)^2 \in \mathbb{C}[z]$ is not simple.
However $J_1(e^{2\pi i t}) = (1,2\pi i) \neq J_1(e^{4\pi it}) = (1,4\pi i).$
\end{example}
\begin{cor}
\label{jetsimple2}
If $P \in K[z]$ is monic, simple, and LCR with degree $n \geq 2,$ then there exists an integer $k \geq 0$ such that for every $s \in \mathbb{T},$ $J_k(\lambda_1),...,J_k(\lambda_n)$ are distinct where $\lambda_1,...,\lambda_n$ are the roots of $P_s \in K_0[z].$
\end{cor}
{\it Proof} Since the roots of $P_s \in K_0[z]$ are analytic functions of $s$ their k-jets are continuous. The result then follows since $\mathbb{T}$ is compact.
\begin{defi}
\label{braid}
For any integer $n \geq 1$ and metric space space $X$ let $C_nX$ denote the metric space consisting of subsets of $X$ having $n$ elements with the Hausdorff metric and let $x \in C_nX.$ The {\it fundamental group} $\pi_1(C_nX,x)$ of $C_n$ with base point $x$ is called the {\it braid group} on $X$ with base point $x$ and denoted by $B_n(X,x).$ Let $S_x$ denote the permutation group on $x.$ We construct a canonical homomorphism $\Phi\, : \, B_n(X,x) \rightarrow S_x$ as follows: let $\phi : [0,1] \rightarrow \mathbb{T}$ be the canonical group homomorphism, let $g \in B_n(X,x),$ and let $h : \mathbb{T} \rightarrow C_nX$ such that $g$ is the homotopy class of $h$ and
$h(0) = h(1) = x.$ Then the composition $f = h \circ \phi \, : \, [0,1] \rightarrow C_nX$ satisfies $f(0) =  f(1) = x$ and therefore induces a continuous map $F \, : \, [0,1] \times X \rightarrow X$ such that
$
    f(t) = \{\, F(t,u) \, : \, u \in x \, \}.
$
We observe that $F(1,\cdot)$ belongs to $S_x$ and does not depend on the representative $h$ for the homotopy class $g.$ Then we define $\Phi(g) = F(1,\cdot).$ The kernel of $\Phi,$ called the {\it pure braid group}, is denoted by $P_n(X,x).$
\end{defi}
We note the well known fact that $C_n\mathbb{C}$ is the configuration space that parameterizes the set of monic polynomials with $n$ distinct roots (\cite{rolfsen}, page 15) and that braid groups arise in both classical mechanics \cite{bas} and quantum physics \cite{ckl}, \cite{kl}. The following result
gives a complete answer to Question \ref{Question 5}.
\begin{theo}
\label{AQ5}
Assume that $P \in K[z]$ is LCR, monic, simple, and $\hbox{deg}(P) = n \geq 2.$ Choose an integer $k \geq 0$ as in Corollary \ref{jetsimple2}, define $h \, : \, \mathbb{T} \rightarrow C_n\mathbb{C}^{k+1}$ by
$$
    h(s) = \{\, J_k(\lambda) \, : \lambda \in K_0, \, P_s(\lambda) = 0 \, \},
$$
let $x = h(0),$ and define $B(P) \in B_n(C_n\mathbb{C}^{k+1},x)$ to be the homotopy class of $h.$
Then $P$ is CR iff $B(P) \in P_n(C_n\mathbb{C}^{k+1},x).$
\end{theo}
{\it Proof} Lemma \ref{jetsimple1} and Corollary \ref{jetsimple2} imply that $B(P)$ is an element in the braid group. The last assertion follows since the roots of $P_s, \, s \in \mathbb{T}$ in $K_0$ can be {\it glued together} to form roots of $P$ in $K$ precisely when $B(P)$ is an element in the pure braid subgroup.
\begin{cor}
\label{corAQ5}
If $h \, : \, \mathbb{T} \rightarrow C_n\mathbb{C}^{k+1}$ is as in Theorem \ref{AQ5} then there there exists a subset $S \subset \mathbb{C}^{k+1}$ of Lebesque measure zero such that the Hermitian product $v*h \, : \, \mathbb{T} \rightarrow C_n\mathbb{C}.$ Therefore $v*h$ gives an element in the braid group $B_n(\mathbb{C},v*x).$
\end{cor}
{\it Proof} The first assertion follows from Sard's theorem and the second assertion is obvious.
\\ \\
Since elements of the braid group are invariant under homotopies of LCR polynomials, it follows that Conjecture \ref{Conj} holds for any
continuous family $P_\kappa$ of LCR polynomials if it holds for any value of $\kappa.$ We note that other homotopy invariants, such as Chern classes, have proved useful in the study of both the Integer Quantum Hall Effect \cite{bes} and the Fractional Quantum Hall Effect \cite{mv}.
\section{Origin of Questions in Quantum Physics}
The questions in this paper arose from a study of the
spectrums of Floquet operators that describe the dynamics of
certain quantum systems. These operators include five families
parameterized by parameters
$\alpha \in (0,1), \lambda > 0, \kappa > 0, \theta \in (0,1).$
\begin{enumerate}
\item Harper (H) operators or Almost Mathieu operators. These self-adjoint operators
can be represented on the Hilbert space $L^2(\mathbb{T})$ with respect to the standard orthonormal basis $\{\xi_n(t) = e^{2\pi i nt}:n \in \mathbb{Z}\}$ by $H(\alpha,\lambda,\theta)\xi_n = \xi_{n-1}+\xi_{n+1} + 2\lambda \cos(2\pi(n\alpha+\theta)).$
\item Unitary Harper (UH) operators $\exp \, [-i \kappa \, H(\alpha,\lambda,\theta) ].$
\item Kicked Harper (KH) operators
$\exp \, \left[-i2\kappa \cos(2\pi t)\right] \, \exp \, \left[-i2\kappa \lambda \cos \left( -i\alpha \frac{d}{dt} + 2\pi \theta \right)\right].$
\item Single Kicked Rotator (SKR) operators
$\exp \, [-i 2 \kappa \cos(2\pi t) ] \, \exp \, \left[\frac{i\alpha}{4\pi} \frac{d^2}{dt^2}\right].$
\item On Resonance Double Kicked Rotator (ORDKR) operators \\
$\exp \, [-i 2 \kappa \cos(2\pi t) ] \, \exp \, \left[-\frac{i\alpha}{4\pi} \frac{d^2}{dt^2}\right]\, \exp \, [-i 2 \kappa \cos(2\pi t) ] \, \exp \, \left[\frac{i\alpha}{4\pi} \frac{d^2}{dt^2}\right].$
\item For rational $\alpha = p/q$ where $p$ and $q$ are coprime integers $\geq 2,$
Mother operators constructed from each of these families of operators by forming
their directed integral over $\theta \in [0,1/q).$
\end{enumerate}
We summarize properties of these Floquet operators that are discussed in detail in \cite{lmwg}.
For $\alpha \in (0,1)$ let $\mathfrak{A}_\alpha$ denote the universal rotation $C^*$-algebra and let $\mathfrak{B}_{\alpha}$ denote the $C^*$-algebra of operators on the Hilbert space $L^2(\mathbb{T})$ generated by (multiplication by) $e^{2\pi i t} \in C(\mathbb{T})$ and the rotation operator $(R_\alpha f)(t) = f(t+\alpha).$ $\mathfrak{B}_\alpha$ is isomorphic to
$\mathfrak{A}_\alpha$ if $\alpha$ is irrational and for rational $\alpha = p/q,$
$\mathfrak{B}_\alpha$ is a homomorphic image of $\mathfrak{A}_\alpha$ and $\mathfrak{B}_\alpha$
is isomorphic to the algebra of matrix valued functions $C(\mathbb{T},\mathbb{C}^{q \times q})$ acting by pointwise multiplication on the Hilbert space $H = L^2(\mathbb{T},\mathbb{C}^q).$ If an operator $F$ corresponds to $M \in C(\mathbb{T},SU(q))$ then
\begin{equation}
\label{speceig}
\hbox{spec}(F) = \bigcup_{t \in \mathbb{T}} \hbox{roots}(zI_q - M(t)),
\end{equation}
where $I_q$ denotes the $q \times q$ identity matrix. Therefore $\hbox{spec}(F)$ consists of the union of at most $q$-disjoint intervals.
\begin{enumerate}
\item For rational $\alpha = p/q,$ the operators 1, 2, 3 and 5 above belong to
$\mathfrak{B}_\alpha$ and their mother operators belong to $\mathfrak{A}_\alpha,$ and the operator 4 belongs to $\mathfrak{B}_{\alpha/2}$ and its mother operators belong to
$\mathfrak{A}_{\alpha/2}.$
\item For rational $\alpha = p/q,$ the spectrum of $H(\alpha,\lambda,\theta))$
and the spectrum of its mother operator equals the union of
$q$ disjoint intervals if $q$ is odd and of $q-1$ disjoint intervals if $q$ is even, (\cite{bs}, Theorems 2 and 3), (\cite{boca01}, Theorem 4.7).
\item The KH operator differs from the UH operator
by $O(\kappa^2)$ so the spectral mapping theorem implies that for fixed rational
$\alpha$ the spectrum of the kicked Harper operator has the same number of bands (disjoint intervals) as described in 2.
\item For irrational $\alpha$ the spectrum of H operator is a Cantor set.
This fact was conjectured by Abzel in 1964 \cite{azbel} and proved by Avila and Jitomirskaya in 2003 \cite{avila}.
\item ORDKR operators were discovered in 2007 by Jiangbin Gong and Jiao Wang
\cite{gongwang07} and further discussed by them and Anders Mouritzen in \cite{wmg}.
They noted that their computed spectrums for rational $\alpha = p/q$ were very close to the spectrums of KH operators with the same parameter values, and that for fixed $\lambda = 1$ and $\kappa$ the Hausdorff distance between their spectrums converged to $0$ as $q$ increased. They also noted that ORDKR systems could be realized using Bose-Einstein condensates, which are described in Appendix B, whereas the KH systems can not be experimentally realized because they require magnetic field strengths five orders of magnitude stronger than the most powerful magnetic fields on Earth (the ones used in MRI devices) and only exist in neutron stars. See (\cite{lmwg}, Appendix A: Physical Considerations and Experimental Realizations) for a detailed discussion of these practical considerations.
\item In \cite{lmwg} these relationships were derived using properties of rotation
$C^*$-algebras. It was explained that KH operators and ORDKR operators are proper homomorphic images of their mother operators if $\alpha$ is rational and that they are isomorphic to their mother operators if $\alpha$ is irrational. Furthermore, their mother operators are unitarily equivalent (under an automorphism in the Brenken-Watatani automorphic representation of the modular group $SL(2,Z)$ acting on $\mathfrak{A}_\alpha)$ and therefore their spectrums are equal. For $\alpha = p/q$ it was proved that the Hausdorff distance between their spectrums approaches $0$ as $q$ increases.
\item In \cite{lmwg} numerical computation of the spectrum of the ORDKR operator for
$\alpha = 2584/4181$ showed that it had a fractal type structure. This supports
our conjecture that the spectrum is a Cantor set if $\alpha$ is irrational.
We briefly discussed approaches to prove our conjecture.
\item In \cite{wg} Cantor type spectrum are observed based on numerical computations
for a class of operators different from the five families mentioned above.
\end{enumerate}
Since for sufficiently small $\kappa$ the multiplicity of roots for polynomial $P_\kappa$ for the ORDKH operator is $\leq 2$ and since the roots of $P_\kappa$ have modulus 1, Corollary \ref{unitary} and Theorem \ref{PPF} imply that
$P_\kappa$ is CR. Therefore if $P_\kappa$ was LCR for all values of $\kappa,$ since the braid group elements are homotopy invariants, Conjecture \ref{Conj} would hold for all $P_\kappa.$ Unfortunately the LCR property is not invariant under ordinary homotopy. Combining Theorem \ref{braid} with modern analytic geometry tools, such as Lojaciewicz's Structure Theorem for Varieties \cite{loj}, (\cite{krantz}, Theorem 5.2.3) and \'{e}tale homotopy, may provide invariants for the $\kappa$ parameterized homotopies to help prove Conjecture \ref{Conj}.
\section{Brief History of Bose-Einstein Condensates}
Bose-Einstein condensates (BEC) arise when a dilute gas of photons or bosonic atoms
are cooled to near absolute zero. Under these conditions a large fraction of
particles occupy the lowest energy state and the gas exhibits {\it weird}
quantum behavior that was first predicted for photons by Satyendra Nath Bose \cite{bose} in 1924 and extended to bosonic matter by Albert Einstein \cite{einstein1, einstein2}. In 1995 Eric Cornell and Carl Weiman \cite{cew} produced the first BEC, consisting of a gas of rubidium atoms cooled to $1.7 \times 10^{-7}\, K,$ at the University of Colorado at Bolder NIST-JILA lab for which they shared the 2001 Nobel Prize in Physics with Wolfgang Ketterle at MIT. Later in 1995 BCE were used to experimentally realize quantum kicked rotators \cite{mrbsr}. In 2010 Jan Klaers, Julian Schmitt, Frank Vewinger and Martin Weitz \cite{ksvw} produced a photon BEC.

\Finishall 

\begin{thebibliography}{10}

\bibitem{abhyankar2} S. S. Abhyankar, {\it Algebraic Geometry for Scientists and Engineers},
American Mathematical Society, Providence, Rhode Island, (1990). [MR 92a:14001]

\bibitem{abhyankar3} S. S. Abhyankar, {\it Irreducibility criteria for germs of analytic
functions of two complex variables}, Advances in Mathematics, 74 (1989), 190--257.

\bibitem{abhyankar4} S. S. Abhyankar, {\it Resolution of Singularities and modular
Galois theory}, Bulletin of the American Mathematical Society, Volume 38, Number 2, (2000), 131--169.

\bibitem{atiyah} M. F. Atiyah and I. G. Macdonald, {\it Introduction to Commutative Algebra},
Addison-Wesley, (1969).

\bibitem{avila} A. Avila and S. Jitomirskaya, {\it The ten martini problem},
Annals of Mathematics, Volume 170, Issue 1, (2009), 302--342. arXiv:math/0503363v1

\bibitem{azbel} M. Y. Azbel, {\it Energy spectrum of a conduction electron in a magnetic
field}, Sov. Phys. JETP 19, (1964), 634.

\bibitem{becker} O. Becker, {\it Eudoxus-Studieren I. Einer voreuklidische Proportionslehre
und ihre Spuren bei Arisoteles und Euklid}, Quellen und Studieren zur Geschichte der
Mathematik, B2, (1933), 211--333.

\bibitem{bezout} E. Bezout, {\it Th\'{e}orie g\'{e}n\'{e}ral des \'{e}quations alg\'{e}briques},
Paris, (1770).

\bibitem{bs} J. Bellisard and B. Simon, {\it Cantor spectrum for the almost Mathieu operator},
Journal of Functional Analysis, 48 (1982), 408--419.

\bibitem{bes} J. Bellisard, A. van Elst, and H. Schultz--Baldes, {\it The noncommutative geometry
of the quantum Hall effect}, Journal of Mathematical Physics, 35 (10), 5373--5451.

\bibitem{bose} S. N. Bose, {\it Plancks Gesetz und Lichtquantenhypothese},
Zeitschrift für Physik, volume 26, number 1, (1924), 178--181.

\bibitem{boca01} F. P. Boca, {\it Rotation $C^*$--Algebras and Almost Mathieu Operators},
Theta, Bucharest (2001).

\bibitem{bas} P. L. Boyland, H. Aref, and M. A. Stremler, {\it Topological fluid mechanics
of stirring}, Journal of Fluid Mechanics, Volume 403, (2000), 277--304.

\bibitem{cew} E. A. Cornell, J. R. Ensher and C. E. Wieman, {\it Experiments in dilute atomic
Bose-Einstein condensation, in Bose-Einstein Condensation in Atomic Gases}, Proceedings of the International School of Physics "Enrico Fermi" Course CXL" (M. Inguscio, S. Stringari and C. E. Wieman, Eds.) Italian Physical Society, (1999), 15-66, cond-mat/9903109.

\bibitem{ckl} G. Chen, L. Kauffman and S. J. Lomonaco, (eds.) {\it The Mathematics of
Quantum Computation and Quantum Technology}, Chapman and Hall/CRC, (2007).

\bibitem{cry} G. Chrystal, {\it Algebra, Vols 1 and 2}, A. and C. Black, Edinburgh, 1889.

\bibitem{eisenbud} D. Eisenbud, {\it Commutative Algebra with a View Toward Algebraic Geometry},
Springer, (1995).

\bibitem{einstein1} A. Einstein, {\it Quantentheorie des einatomigen idealen Gases},
Sitzungsberichte der Preussischen Akademie der Wissenschaften, (1924), 261--267.

\bibitem{einstein2} A. Einstein, {\it Quantentheorie des einatomigen idealen Gases. 2.
Abhandlung}, Sitzungsberichte der Preussischen Akademie der Wissenschaften, (1925), 3--14.

\bibitem{edwards} H. M. Edwards, {\it Essay 4.4 Newton's Polygon,
pages 132--141 in Essays in Constructive Mathematics},
Springer, New York, (2005).

\bibitem{gongwang07} J. Gong and J. Wang,
{\it Quantum diffusion dynamics in nonlinear systems: a modified kicked-rotor model}, Phys. Rev. E, Volume 76, Issue 3, 036217, 5 pages, (2007).

\bibitem{gunning} R. C. Gunning and H. Rossi, {\it Analytic Functions of Several
Complex Variables}, Prentice Hall, (1965).

\bibitem{heath} T. H. Heath, {\it The Thirteen Books of Euclid's Elements}, 2nd ed., Dover,
(1956).

\bibitem{hensel} K. Hensel, {\it Theorie der algebraischen zahlen}, Teubner, Leipzig, (1908).

\bibitem{hironaka} H. Hironaka, {\it Resolution of singularities of an algebraic variety
over a field of characteristic zero, I, II}, Annals of Mathematics, (2) 79 (1964), 109--326.

\bibitem{kl} L. H. Kauffman and S. J. Lomonaco, {\it Unitary Representations of the
Artin Braid Groups and Quantum Algorithms for Colored Jones Polynomials and the Witten-Reshetikhin-Turaev Invariants, in Topology and Physics - Proceedings of the Nankai International Conference in Memory of Xia-Song Lin, ed. by K. Lin, Z. Wang, and W. Shang}, Nankai Tracts in Mathematics, Vol. 12, (2008), 172-194.

\bibitem{krantz} S. G. Krantz and H. R. Parks, {\it A Primer of Real Analytic Functions},
Birkh\"{a}ser, (1992).

\bibitem{ksvw} J. Klaers, J. Schmitt, F. Vewinger and M. Weitz,
{\it Bose-Einstein condensation of photons in an optical microcavity},
Nature 468 (2010), 545-548.

\bibitem{kunz} E. Kunz, {\it Introduction to Commutative Algebra and Algebraic Geometry},
Birkh\"{a}user, (1985).

\bibitem{lang} S. Lang, {\it Algebra}, 2nd ed., Addison Wesley, (1984).

\bibitem{lmwg} W. Lawton, A. S. Mouritzen, J. Wang, J. B. Gong,
{\it Spectral relationships between kicked Harper and on-resonance double kicked rotor operators}, Journal of Mathematical Physics, Volume 50, 032103, 26 pages, (2009).
arXiv:0807.4276

\bibitem{loj} S. Lojaciewicz, {\it Complex Analytic Geometry}, Birkh\"{a}ser, (1991).

\bibitem{mv} M. Marcolli and M. Varghese, {\it Towards the fractional quantum Hall effect:
a noncommutative geometry perspective, pages 235-261 in Noncommutative Geometry and Number Theory, Editors: C. Consani and M. Marcolli}, Aspects of Mathematics, Vieweg Verlag, Wiesbaden, (2006). cond-mat/0502356.

\bibitem{mrbsr} F. L. Moore, J. C. Robinson, C. F. Bharucha, B. Sundaram and M. G.
Raizen, {\it Atom optics realization of the quantum delta-kicked rotator}, Phys. Rev. Lett. 75, (1995), 4598--4601.

\bibitem{newt} I. Newton, {\it Geometria Analitica}, (1660).

\bibitem{puis} V. A. Puiseux, {\it Recherches sur les fonctions alg\'{e}briques},
J. Math. Pures Appl. 15(1850) 365--480.

\bibitem{rolfsen} D. Rolfsen, {\it Tutorial on the braid groups, pages 1-31 in
Braids: introductory lectures on braids, configurations and their applications, edited by
A. J. Berrick, F. R. Cohen, and E. Hanbury, Lecture Note Series, Volume 19,
Institute of Mathematical Sciences, National University of Singapore}, World Scientific, (2010).

\bibitem{sylvester} J. J. Sylvester, {\it On a general method of determining by mere inspection
the derivations from two equations of any degree}, Philosophical Magazine 16, (1840), 132--135.

\bibitem{wmg} J. Wang, A. S. Mouritzen, J. Gong, {\it Quantum control of ultra-cold
atoms: uncovering a novel connection between two paradigms of quantum nonlinear dynamics}, Journal of Modern Optics, Volume 56, Issue 6, (2008), 722--728.

\bibitem{wg} J. Wang and J. Gong, {\it Generating a fractal butterfly Floquet spectrum
in a class of driven SU(2) systems},
Physical Review E, Volume 81, Issue 2, (2010), 026204.

\end{thebibliography}
\end{document}